\documentclass{agtart_a}
\pdfoutput=1

\usepackage[all]{xy}
\usepackage{pinlabel}

\title{A categorification for the Tutte polynomial}

\author{Edna F Jasso-Hernandez}
\givenname{Edna F}
\surname{Jasso-Hernandez}
\address{Department of Mathematics\\
George Washington University\\\newline
Washington, DC 20052\\USA}
\email{fanny@gwu.edu}
\urladdr{}

\author{Yongwu Rong}
\givenname{Yongwu}
\surname{Rong}
\email{rong@gwu.edu}
\urladdr{}

\volumenumber{6}
\issuenumber{}
\publicationyear{2006}
\papernumber{71}
\startpage{2031}
\endpage{2049}

\doi{}
\MR{}
\Zbl{}

\keyword{Khovanov homology}
\keyword{Tutte polynomial}
\keyword{categorification}
\keyword{graph polyomial}
\subject{primary}{msc2000}{05C15}
\subject{secondary}{msc2000}{57M27}
\subject{secondary}{msc2000}{55N35}

\received{24 January 2006}
\revised{}
\accepted{5 July 2006}
\published{19 November 2006}
\publishedonline{19 November 2006}
\proposed{}
\seconded{}
\corresponding{}
\editor{CPR}
\version{}

\arxivreference{math.CO/0512613}

\let\xysavmatrix\xymatrix
\def\xymatrix{\disablesubscriptcorrection\xysavmatrix}
\AtBeginDocument{\let\bar\wbar\let\tilde\wtilde\let\hat\what}

\makeatletter
\def\cnewtheorem#1[#2]#3{\newtheorem{#1}{#3}[section]
\expandafter\let\csname c@#1\endcsname\c@theorem}
\makeatother

\newtheorem{theorem}{Theorem}[section]
\cnewtheorem{lemma}[theorem]{Lemma}
\cnewtheorem{proposition}[theorem]{Proposition}
\cnewtheorem{corollary}[theorem]{Corollary}
\theoremstyle{remark}
\cnewtheorem{remark}[theorem]{Remark}

\cnewtheorem{acknowledgement}[theorem]{Acknowledgement}
\newtheorem{exa}{Example}

\begin{document}

\begin{abstract} 
For each graph, we construct a bigraded chain complex whose graded
Euler characteristic is a version of the Tutte polynomial. This work
is motivated by earlier work of Khovanov, Helme-Guizon and Rong, and
others.
\end{abstract}

\maketitle

\section{Introduction}\label{sec1}

 In \cite{K00}, Khovanov introduced a graded homology theory for
classical links and showed it yields the Jones polynomial by taking
the graded Euler characteristic. This construction has sparked a good
deal of interests in recent years. In \cite{HR04}, Helme-Guizon and
Rong constructed a graded homology theory for graphs.  The graded
Euler characteristic of the homology groups is the chromatic
polynomial of the graph.

It is natural to ask if similar constructions can be made for
other graph polynomials, especially the Tutte polynomial, which is
universal among graph invariants satisfying the
deletion-contraction rule. In this paper, we give such
constructions for the Tutte polynomial. More precisely, for each
graph $G$ we define bigraded homology groups whose Euler
characteristic is a variant of the Tutte polynomial. Our
construction is different from the one by Khovanov and Rozansky
for categorification for the Homflypt polynomial
\cite{KR04,KR05}.

Our construction starts with rewriting the Tutt polynomial using a
state sum which is more amenable for a chain complex set up. This is
done in \fullref{sec2}.
The chain groups will be built on two basic
algebraic objects: a bigraded algebra $A$, and a bigraded
$\mathbb{Z}$--module $B$. The differential will depend on the
multiplication on $A$ and a sequence of graded homomorphisms $f_k\co 
B^{\otimes k} \rightarrow B^{\otimes k+1}$ ($k=0, 1, 2, \cdots$).  In
\fullref{sec3} we construct our homology groups using an obvious choice of
$A, B$ and $f_k$.  A more general construction is shown in Section
4. These homology groups satisfy a long exact sequence which we
explain in \fullref{sec5}. In \fullref{sec6}, we discuss additional properties,
including a functorial property.  Some computational examples are
given in \fullref{sec7}.

We wish to thank Laure Helme-Guizon for her helpful discussions. The
second author is partially supported by NSF grant DMS-0513918.

\section{The Tutte polynomial}\label{sec2}
We recall some basic properties of the Tutte polynomial. Let $G$ be
a graph with vertex set $V(G)$ and edge set $E(G)$. Given an edge
$e\in E(G)$, let $G-e$ denote the graph obtained from $G$ by
deleting the edge $e$, let $G/e$ denote the graph obtained by
contracting $e$ to a vertex. Recall that $e$ is called a loop if $e$
joins a vertex to itself, $e$ is called an isthmus if deleting $e$
from $G$ increases the number of components of the graph. The Tutte
polynomial of $G$, denoted by $T(G; x, y)$, is uniquely defined by
the following axioms: 

\begin{itemize}
\item $ T(G; x, y)=T(G-e; x, y)+ T(G/e; x, y)$, if $e$ is not a loop or
isthmus.

\item $T(G; x, y)=xT(G-e; x, y)$, if $e$ is an isthmus.

\item $T(G; x, y)=y T(G/e; x, y), \mbox{ \ if $e$ is a loop. }
\label{loop}$

\item $T(G;x,y)=1$ if $G$ has no edges. 
\end{itemize}

 It is obvious that $T(G; x,y)$ is a 2--variable
polynomial in $x$ and $y$. Furthermore, it has a closed form
described below. First we introduce some notations. Let $s\subset
E(G)$. The rank of $s$, denoted by $r(s)$, is defined by
$r(s)=|V(G)|-k(s) $ where $k(s)$ is the number of connected
components of the graph $[G:s]$ having vertex set $V(G)$ and edge
set $s$. We have the following well-known state sum formula (see,
eg Welsh \cite{W}):
\[ T(G; x,y)=\sum_{s \subset E(G)} (x-1)^{r(E)-r(s)} (y-1)^{|s|-r(s)}
\]
For our construction, we need to rewrite $T(G; x, y)$ in a new form
that is easier to work with. Since $x-1=-(1-x)$, $y-1=-(1-y)$, we
have
\begin{eqnarray*}
T(G;x,y)&=&\sum_{s \subset E(G)}
(-1)^{r(E)-r(s)}(1-x)^{r(E)-r(s)}(-1)^{|s|-r(s)}(1-y)^{|s|-r(s)}\\
&=&(-1)^{r(E)} \sum_{s\subset E(G)}(-1)^{|s|}(1-x)^{r(E)-r(s)}
(1-y)^{|s|-r(s)}.
\end{eqnarray*}
But $r(E)-r(s)=|V(G)|-k(E) - (|V(G)|-k(s))=-k(E)+k(s)$. Note that
$k(s)$ is the number of components of the graph $[G:s]$ which is the
the $0^{\mbox{\tiny th}}$ Betti number of the (underlying space of
the) graph $[G:s]$. We denote $k(s)$ by $b_0(s)$, and hence
$r(E)-r(s)=-k(E)+b_0(s)$. For $|s|-r(s)$, we have
$|s|-r(s)=|s|-(|V(G)|-k(s))=- (|V(G)|-|s|)+k(s)$. But $|V(G)|-|s|$
is the Euler characteristic of the graph $[G:s]$. Let $b_1(s)$
denote the first betti number of $[G:s]$. We have $|V(G)|-|s|=$\break $\chi
([G:s])=b_0(s)-b_1(s)$, and $k(s)=b_0(s)$. It follows that
$|s|-r(s)=b_1(s)$. Therefore
\begin{eqnarray*}
T(G;x,y)&=&(-1)^{r(E)}\sum_{s\subset
E(G)}(-1)^{|s|}(1-x)^{-b_0(E)+b_0(s)} (1-y)^{b_1(s)}\\
&=&(-1)^{r(E)}(1-x)^{-b_0(E)}\sum_{s\subset E(G)} (-1)^{|s|}
(1-x)^{b_0(s)} (1-y)^{b_1(s)}.
\end{eqnarray*}
We have proved:

\begin{proposition} \label{T-hat}
$T(G;x,y)=(-1)^{r(E)}(1-x)^{-b_0(E)}\widehat{T}(G;-x,-y)$, where\break
$\widehat{T}(G;x,y)=\sum_{s\subset E(G)} (-1)^{|s|} (1+x)^{b_0(s)}
(1+y)^{b_1(s)}$, and $b_i(s)=b_i([G:s])$ is the $i^{\mbox{\tiny
th}}$ Betti number of $[G:s]$ ($i=0, 1$).
\end{proposition}

It is the 2--variable polynomial $\widehat{T}(G;x,y)$ that we will
categorify. This does recover the Tutte polynomial because of the
following lemma.

\begin{lemma}
The Tutte polynomial $T(G)$ is determined by $\widehat{T}(G)$.
\end{lemma}

\begin{proof} We make the change of variables $u=1+x$, $v=1+y$ which
turns $\widehat{T}(G;x,y)$ into $\sum_{s\subset E(G)} (-1)^{|s|}
u^{b_0(s)} v^{b_1(s)}$, which we denote by $\widetilde{T}(G;u,v)$.
Each subset $s\subset E(G)$ yields a term in $\widetilde{T}(G;u,v)$.
Define the complexity of each such term to be $comp(s)=(-b_0(s),
b_1(s))$ with the dictionary order. By \fullref{adding edge}
below, the complexity goes up when adding an edge to $s$. It follows
that when $s=\emptyset$, we have the minimum term which is
$u^{|V(G)|}$, and when $s=E(G)$, we have the maximum term which is
$(-1)^{|E|}u^{b_0(E)} v^{b_1(E)}$. Therefore, we can recover
$b_0(E)$, and $r(E)$ (which is $|V(G)|-b_0(E)$) from
$\widehat{T}(G;x,y)$. By \fullref{T-hat}, $T(G;x,y)$ can be
recovered from $\widehat{T}(G;x,y)$.
\end{proof}

\begin{lemma} \label{adding edge}
Let $s$ be a subset of $E(G)$, and $e$ be an edge not in $s$. Then
one of the following two cases occurs.
\begin{itemize}
\item[{\rm(i)}] If $e$ joins two components of $[G:s]$, then $b_0(s
\cup \{ e\})=b_0(s)-1$ and $b_1(s) = b_1(s \cup \{ e\})$. 

\item[{\rm(ii)}] If $e$ connects a component of $[G:s]$ to itself, then  $b_0(s)
= b_0(s \cup \{ e\})$ and $b_1(s \cup \{ e\})=b_1(s)+1$.
\end{itemize}
\end{lemma}

This lemma follows from standard algebraic topology and its proof is
omitted.

\section{The chain complex}\label{sec3}

\subsection{Algebra background}
 Let $A$ be a commutative algebra over a commutative ring. For our purpose,
 the ring will be $\Z$.
Recall that $A$ is called  a graded algebra if it can be written as
a direct sum $A=\oplus_{i\in \Z}A_i$ where each $A_i$ is closed
under addition, and $A_i A_j \subset A_{i+j}$ (ie, $a_i a_j\in A_{i+j}$
for all $a_i\in A_i, a_j\in A_j$). The elements in $A_i$ are called
homogeneous elements of degree $i$. Thus the condition $A_i A_j
\subset A_{i+j}$ is equivalent to the condition that the degree is
additive under multiplication.

The same definition can be made for modules over a ring by simply
dropping the additivity condition of degree. Thus a $\Z$--module $M$
is a graded module if we write it as a direct sum of submodules
$M=\oplus_{i\in \Z}M_i$ where elements of $M_i$ are called
homogeneous elements of degree $i$. By definition, a graded algebra
is automatically a graded module over the same ring.

Let $M=\oplus_{i\in \Z}M_i$ be a graded $\mathbb{Z}$--module. The
{\it graded dimension} of $M$ is the power series
\[
q\dim M:=\sum_i q^i \cdot \dim_{\mathbb{Q}} (M_i \otimes \mathbb{Q})
\]
The same definition holds for a graded algebra.

An obvious generalization for graded algebras (modules) can be made
by allowing the grading index $i$ lying in an arbitrary abelian
group. In this paper, we are interested in the case when the group
is $\Z \oplus \Z$. Such an algebra (resp. module) will be called a
$\Z \oplus \Z$ --graded, or a bigraded algebra (resp. module). By
definition, a bigraded algebra is an algebra $A$ with a
decomposition $A=\oplus_{(i,j)\in \Z \oplus \Z}A_{i,j}$ where each
$A_{i,j}$ is closed under addition and $A_{i_1, j_1} A_{i_2, j_2}
\subset A_{i_1+i_2, j_1+j_2}$ for all $(i_1,j_1), (i_2,j_2) \in \Z
\oplus \Z$. A bigraded module is a module $M$ with a decomposition
$M=\oplus_{(i,j)\in \Z \oplus \Z }M_{i,j}$ where elements of
$M_{i,j}$ are called homogeneous elements with degree $(i,j)$. The
graded dimension of $M$ is the 2--variable power series
\[
q\dim M:=\sum_{i,j} x^i y^j \cdot \dim_{\mathbb{Q}} (M_{i,j} \otimes
\mathbb{Q})
\]

\subsection{The construction -- a specific one}\label{sec3.2}

Let $A=\mathbb{Z}[x]/(x^2), B=\Z[y]/(y^2)$ where $\deg x=(1,0), \deg
y=(0,1)$. Then $A$ and $B$ become bigraded algebras with $ q\dim
A=1+x, q\dim B=1+y$. Note that $A^{\otimes m}\otimes B^{\otimes n}$
is a bigraded $\mathbb{Z}$--module whose graded dimension is $q\dim
A^{\otimes m}\otimes B^{\otimes n}=(1+x)^m(1+y)^n$. We are not going
to use the algebra structure on $B$, but its graded module structure
will be needed.

Note that the letter $x$ has two different meanings. First, it is an
element in the algebra $A$. Second, it is the variable in the power
series $q\dim A$. The same is true for the letter $y$. We believe
such abuse of notation is convenient and will not lead to confusion.

Now, let $G$ be a graph with $|E(G)|=n$. We fix an ordering on
$E(G)$ and denote the edges by $e_1, \cdots, e_n$. Consider the
$n$--dimensional cube $\{ 0, 1\}^E = \{ 0, 1\}^n$. Each vertex
$\alpha$ of this cube corresponds to a subset $s=s_{\alpha}$ of $E$,
where $e_i \in s_{\alpha}$ if and only if $\alpha_i=1$. The height
$|\alpha|$ of $\alpha$, is defined by $|\alpha|=\sum \alpha_i$,
which is also equal to the number of edges in $s_{\alpha}$.

For each vertex $\alpha$ of the cube, we associate the graded
$\mathbb{Z}$--module $C^{\alpha}(G)$ as follows. Consider $[G:s]$,
the graph with vertex set $V(G)$ and edge set $s$. We assign a copy
of $A$ to each component of $[G:s]$ and then take tensor product
over the components. Let $A^{\alpha}(G)$ be the resulting graded
$\mathbb{Z}$--module, with the induced grading from $A$.  Therefore,
$A^{\alpha}(G) \cong A^{\otimes k}$ where $k=b_0([G:s])$ is the
number of components of $[G:s]$. Next, let $B^{\alpha}(G)=B^{\otimes
l}$ where $l=b_1([G:s])$ is the first Betti number of $[G:s]$ (note
that there is no specific order on the tensor factors here). We
define $C^{\alpha}(G)=A^{\alpha}(G)\otimes B^{\alpha}(G)$. Then we
define the $i^{\mbox{\tiny th}}$ chain group
$C^i(G):=\oplus_{|\alpha|=i} C^{\alpha}(G)$. This defines the chain
groups of our complex.

As a notational remark, we can also denote $C^{\alpha}(G)$ by
$C^{s}(G)$, because of the one-to-one correspondence between $s$ and
$\alpha$.

Next, we define the differential maps $d^i\co  C^i(G)\rightarrow
C^{i-1}(G)$. We need to make use of the edges of the cube $\{ 0,
1\}^E$. Each edge $\xi$ of $\{ 0, 1\}^E$ can be labeled by a
sequence in $\{ 0, 1,
*\}^E$ with exactly one $*$. The tail of the edge is obtained by
setting $*=0$ and the head is obtained by setting $*=1$. The height
$|\xi|$ is defined to the height of its tail, which  is also equal
to the number of 1's in $\xi$.

Given an edge $\xi$ of the cube, let $\alpha_1$ be its tail and
$\alpha_2$ be its head. Let $e$ be the corresponding edge in $G$ so
that $s_2=s_1\cup \{ e\}$. We define the per-edge map $d_{\xi}\co 
C^{\alpha_1}(G)\rightarrow C^{\alpha_2}(G)$ based on the two cases
in \fullref{adding edge}.

{\bf Case 1}\qua  $e$ joins a component of $[G:s_1]$ to itself. The
components of $[G:s_2]$ and $[G:s_1]$ naturally correspond to each
other, and therefore $A^{\alpha_1}(G)= A^{\alpha_2}(G)$. Let
$d_{\xi}^A\co  A^{\alpha_1}(G)\rightarrow A^{\alpha_2}(G)$ be the
identity map. We also have $B^{\alpha_2}=B^{\alpha_1}\otimes B$. Let
$d_{\xi}^B\co  B^{\alpha_1}(G)\rightarrow B^{\alpha_2}(G)$ be the
homomorphism sending $b\in B^{\alpha_1}(G)$ to $b\otimes 1\in
B^{\alpha_1}(G)\otimes B=B^{\alpha_2}(G)$. The per-edge map
$d_{\xi}\co  C^{\alpha_1}(G)\rightarrow C^{\alpha_2}(G)$ is defined by
$d_{\xi}=d_{\xi}^A\otimes d_{\xi}^B\co  A^{\alpha_1}(G)\otimes
B^{\alpha_1}(G) \rightarrow A^{\alpha_2}(G)\otimes B^{\alpha_2}(G)$.

{\bf Case 2}\qua $e$ joins two different components of $[G:s_1]$.
 Let $E_1, E_2, \cdots, E_{k}$ be the
components of $[G:s_1]$, where $E_1$ and $E_2$ are connected by $e$.
Then the components of $[G:s_2]$ are $E_1\cup E_2 \cup \{ e\}, E_3,
\cdots, E_{k}$. We define $d_{\xi}^A$ to be the identity map on the
tensor factors coming from $E_3, \cdots, E_{k}$, and $d_{\xi}^A$ on
the tensor factors coming from $E_1, E_2$ to be the multiplication
map $A \otimes A \rightarrow A$ sending $a_1\otimes a_2$ to $a_1
a_2$. For $d_{\xi}^B$, we have $B^{\alpha_2}=B^{\alpha_1}$ and we
define $d_{\xi}^B$ to be the identity map. Again, the per-edge map
is defined by $d_{\xi}=d_{\xi}^A\otimes d_{\xi}^B$.

Now, we define the differential $d^i\co C^i(G) \rightarrow C^{i+1}(G)$
by $d^i = \sum_{|\xi|=i} (-1)^{\xi} d_{\xi}$, where $(-1)^{\xi}
=(-1)^{\sum_{i<j} \xi_i} $ and $j$ is the position of $*$ in $\xi$.

\begin{figure}[ht!]
\vspace{2mm}
\cl{
\labellist\tiny
\pinlabel 000 [tr] at 191 398
\pinlabel 100 [tr] at 457 566
\pinlabel 010 [tr] at 457 398
\pinlabel 001 [tr] at 457 225
\pinlabel 110 [tr] at 733 566
\pinlabel 101 [tr] at 733 398
\pinlabel 011 [tr] at 733 225
\pinlabel 111 [tr] at 1001 398
\pinlabel $i=0$ [b] at 102 571
\pinlabel $i=1$ [b] at 365 571
\pinlabel $i=2$ [b] at 640 571
\pinlabel $i=3$ [b] at 906 571
\pinlabel 0 [r] at 10 8
\pinlabel $C^0(G)$  at 89 8
\pinlabel $d^0$  [b] at 218 8
\pinlabel $C^1(G)$  at 363 8
\pinlabel $d^1$  [b] at 489 8
\pinlabel $C^2(G)$  at 633 8
\pinlabel $d^2$  [b] at 762 8
\pinlabel $C^3(G)$  at 902 8
\pinlabel 0 at 991 8
\pinlabel $A{\otimes}A$ [b] at 89 276
\pinlabel $A$ [b] at 363 443
\pinlabel $A$ [b] at 363 276
\pinlabel $A{\otimes}A{\otimes}B$ [b] at 363 104
\pinlabel $A{\otimes}B$ [b] at 633 443
\pinlabel $A{\otimes}B$ [b] at 633 276
\pinlabel $A{\otimes}B$ [b] at 633 104
\pinlabel $A{\otimes}B{\otimes}B$ [b] at 902 276
\pinlabel $\oplus$ at 363 421
\pinlabel $\oplus$ at 363 250
\pinlabel $\oplus$ at 633 421
\pinlabel $\oplus$ at 633 250
\pinlabel* $d_{*00}$ [tl] at 224 428
\pinlabel* $d_{0*0}$ <0pt,1pt> [b] at 223 338
\hair 1pt
\pinlabel $d_{00*}$ [bl] at 224 239
\pinlabel* $d_{1*0}$ <0pt,2pt> [b] at 496 509
\pinlabel* $d_{*10}$ <3pt,3pt> [br] at 522 458
\pinlabel* $d_{10*}$ <4pt,0pt> [tr] at 518 365
\pinlabel* $d_{*01}$ <6pt,8pt> [r] at 515 299
\pinlabel* $d_{01*}$ <3pt,-1pt> [r] at 519 203
\pinlabel* $d_{0*1}$ [t] at 497 150
\pinlabel* $d_{*11}$ [tl] at 771 235
\pinlabel $d_{11*}$ [bl] at 769 446
\pinlabel $d_{1*1}$ [b] at 768 343
\endlabellist
\psfig{figure=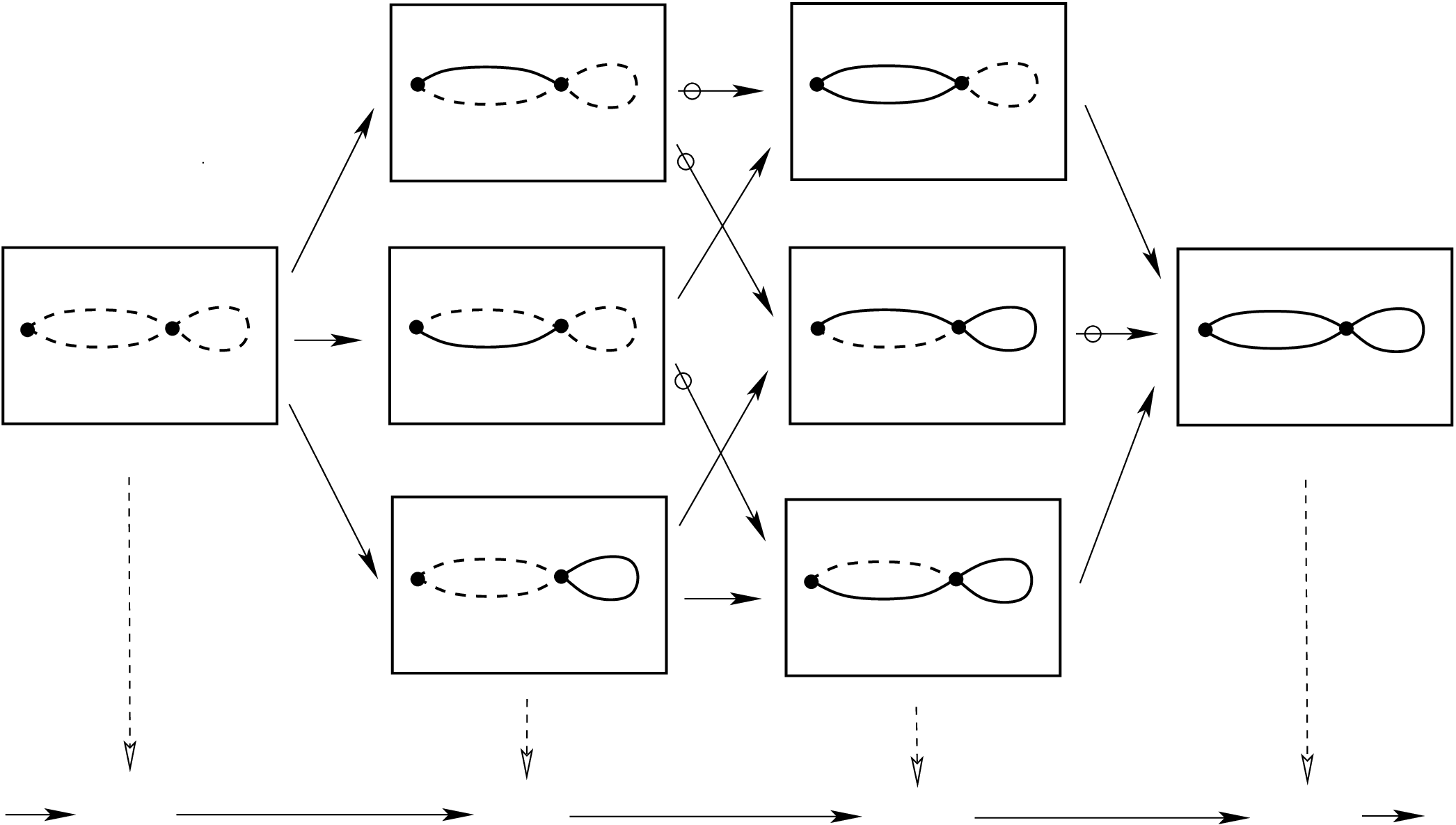,width=.99\hsize}}
\caption{}\label{fig1}
\end{figure}

To illustrate our construction consider the graph with two vertices,
two parallel edges and a loop attached to one of its vertices. Let
us label the edges of $G$ as follows:
\cl{
\labellist\small
\pinlabel $G=$ [r] at 37 45
\pinlabel $e_1$ [b] at 85 61
\pinlabel $e_2$ [t] at 85 26
\pinlabel $e_3$ [t] at 170 32
\endlabellist
\psfig{figure=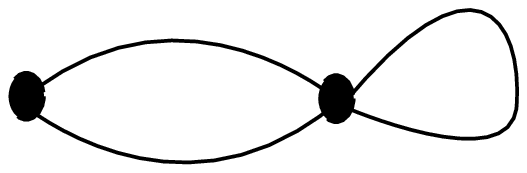,height=1.5 cm}}
The representation of the 3--dimensional cube and the chain complex
is given in \mbox{\fullref{fig1}}.
In each rectangular box, the right upper corner has the
sequence corresponding to $\alpha$, the center has the graph $[G:s]$
and the algebra $C^\alpha (G)$ is at the bottom. Taking direct sum on
each column gives the chain group $C^i(G)$ in the bottom row.

\newpage
In \fullref{fig1} we have represented the per-edge maps $d_ \xi$. The
arrows with a circle represent the maps for which $(-1)^{\xi}= -1$.
\vspace{2pt}

For this particular example we have that both $d_{*00}$ and
$d_{0*0}$ map $a_1 \otimes a_2 \mapsto a_1 a_2$,  $d_{00*}$ maps
$a_1 \otimes a_2 \mapsto a_1 \otimes a_2 \otimes 1_B $. Also
$d_{1*0}$, $d_{10*}$ and $d_{01*}$ map $a \mapsto - a \otimes 1_B$;
 $d_{*10}$ maps $a \mapsto a \otimes 1_B$; both $d_{0*1}$ and $d_{*01}$ map
  $a_1 \otimes a_2 \otimes
b \mapsto a_1 a_2 \otimes b$. Finally $d_{11*}$ and $d_{*11}$ map $a
\otimes b \mapsto a \otimes b \otimes 1_B$, $d_{1*1}$ maps $a \otimes b \mapsto -a \otimes b \otimes 1_B$.
\vspace{2pt}

\begin{theorem}\label{chain complex}$\phantom{99}$

{\rm(a)}\qua $0 \rightarrow C^0(G)
\stackrel{d^0}{\rightarrow}C^1(G) \stackrel{d^1}{\rightarrow} \cdots
\stackrel{d^{n-1}}{\rightarrow} C^n(G) \rightarrow 0$ is a chain
complex of  bigraded modules whose differential
is degree preserving. Denote this chain complex by $C(G)$. 

{\rm(b)}\qua The cohomology groups $H^i(G)$ are independent of the ordering
of
the edges of $G$, and therefore are invariants of the graph $G$.
In fact, the isomorphism class of $C(G)$ is an invariant of $G$.

{\rm(c)}\qua $\chi_{q}(C(G))=\mathrel{\mathop{\sum }\limits_{0\leq i\leq
n}}(-1)^{i} q\dim (H^{i})=\mathrel{\mathop{\sum }\limits_{0\leq
i\leq n}} (-1)^{i} q\dim (C^{i})=\widehat{T}(G;x,y)$
\end{theorem}
\vspace{2pt}

\begin{proof}
(a)\qua The map $d$ is obviously linear. It is also degree preserving
since it is built on two basic maps: the multiplication on $A$ and
the map $b \rightarrow b\otimes 1$, both being degree preserving.
\vspace{2pt}

It remains to show that $d^2=0$. Let $s\subset E(G)$. Consider the
result of adding two edges $e_k$ and $e_j$ to $s$ where $k<j$. It is
enough to show that the following diagram commutes.
\begin{equation}  \label{diag1}
\begin{array}{ccccc}
&  & C^{ s\cup \left\{ e_k\right\} }(G)  &  &  \\
& \begin{picture}(10,10)(0,0)
  \linethickness{1pt}
  \put(-3,-3){\vector(1,1){18}}
  \put(-35,9){$d_{...*...0...}$}
  \end{picture}
   &  & \begin{picture}(10,10)(0,0)
        \linethickness{1pt}
        \put(-8,15){\vector(1,-1){18}}
        \put(3,10){$d_{...1...*...}$}
        \end{picture}
     &  \\
C^{S}(G)  &  &  &  & C^{s \cup \left\{ e_k, e_j \right\}}(G) \\
& \begin{picture}(10,10)(0,0)
  \linethickness{1pt}
  \put(-3,8){\vector(1,-1){18}}
  \put(-35,-4){$d_{...0...*...}$}
  \end{picture}
 &  & \begin{picture}(10,10)(0,0)
      \linethickness{1pt}
      \put(-8,-5){\vector(1,1){18}}
      \put(2,-5){$d_{...*...1...}$}
      \end{picture}
    &  \\
&  & C^{ s\cup \left\{ e_j\right\} }(G)  &  &
\end{array}
\end{equation}
where $d_{...*...0...}$ means that the vector $\xi$
has a star in the $k$-th position and a zero in the $j$-th position
and the same convention applies to the remaining maps.
 We observe that in  Diagram \ref{diag1} all the per-edge maps will have the same
sign except for exactly one. More precisely $d_{...1...*...}$ will
have a different sign if the number of $1$'s between the star and
the zero is even, and $d_{...0...*...}$ will differ in sign if the
number of 1's is odd. This, along with the commutativity of Diagram
\ref{diag1}, implies $d^{i+1} d^i=0$.
\vspace{2pt}

The commutativity of Diagram \ref{diag1} follows from some tedious
but straight forward checking, based on various ways $e_k$ and $e_j$
join the components of $[G:s]$. Let $e_k$ connect $E_k$ to $F_k$,
and $e_j$ connect $E_j$ to $F_j$, where $E_k, E_j, F_k, F_j$ are
(not necessarily distinct) components of $[G:s]$. Consider the set
$C=\{ E_k, F_k, E_j, F_j\}$. Up to symmetry (ie, interchange of
$k$ and $j$,  and of $E$ and $F$), there are seven possibilities
shown below.

\begin{enumerate}
\item $|C|=1$. We have $E_k=F_k=E_j=F_j$.

\item  $|C|=2$, and $E_k=F_k$, $E_j=F_j$.

\item  $|C|=2$, and $E_k=F_k=E_j$.

\item  $|C|=2$, and $E_k=E_j$, $F_k=F_j$.

\item  $|C|=3$, and $E_k=F_k$.

\item  $|C|=3$, and $F_k=E_j$.

\item  $|C|=4$.
\end{enumerate}

As an example, we check the commutativity for Case 3. This is
depicted in the following diagram.
\begin{equation*}
\begin{array}{ccccccccccc}

&  & A \otimes A \otimes B &  & &
 \hspace{.5cm} &&  &  a_1 \otimes a_2 \otimes 1_B &  &  \\

 & \begin{picture}(10,10)(0,0)
  \linethickness{1pt}
  \put(-3,-3){\vector(1,1){18}}
  \put(-35,9){$d_{...*...0...}$}
 \end{picture}
   &  & \begin{picture}(10,10)(0,0)
        \linethickness{1pt}
        \put(-8,15){\vector(1,-1){18}}
        \put(3,10){$d_{...1...*...}$}
        \end{picture}
     &  & && \begin{picture}(10,10)(0,0)
  \linethickness{1pt}
  \put(-3,-3){\vector(1,1){18}}
  \put(-35,9){$d_{...*...0...}$}
 \end{picture}
   &  & \begin{picture}(10,10)(0,0)
        \linethickness{1pt}
        \put(-8,15){\vector(1,-1){18}}
        \put(3,10){$d_{...1...*...}$}
        \end{picture}
     & \\

\ A \otimes A &  &  &  & A \otimes B
 & \hspace{-.5cm} & a_1 \otimes a_2&  &  &  &  a_1 a_2 \otimes 1_B \\

& \begin{picture}(10,10)(0,0)
  \linethickness{1pt}
  \put(-3,8){\vector(1,-1){18}}
  \put(-35,-4){$d_{...0...*...}$}
  \end{picture}
 &  & \begin{picture}(10,10)(0,0)
      \linethickness{1pt}
      \put(-8,-5){\vector(1,1){18}}
      \put(2,-5){$d_{...1...*...}$}
      \end{picture}
    & & & & \begin{picture}(10,10)(0,0)
  \linethickness{1pt}
  \put(-3,8){\vector(1,-1){18}}
  \put(-35,-4){$d_{...0...*...}$}
  \end{picture}
 &  & \begin{picture}(10,10)(0,0)
      \linethickness{1pt}
      \put(-8,-5){\vector(1,1){18}}
      \put(2,-5){$d_{...1...*...}$}
      \end{picture}
    &  \\

&  &  A &  & & &  & &  a_1  a_2 &  & \\
\end{array}
\end{equation*}

(b)\qua The proof is similar to \cite[Theorem 12]{HR04}. Each
permutation of the edges of $G$ is a product of transpositions of
the form $(k, k+1)$. An explicit isomorphism can be constructed
for each such transposition. In fact, this shows that the
isomorphism class of the chain complex is an invariant of the
graph.

(c)\qua First, a standard homological algebra argument
 shows that $$\sum_{0\leq i\leq n}(-1)^{i} q\dim
(H^{i}(G))=\sum_{0\leq i\leq n} (-1)^{i} q\dim (C^{i}(G)).$$ Next,
each $C^i(G)$ is a direct sum of $C^{\alpha}(G)$ where $\alpha$
corresponds to $s\subset E(G)$ with $|s|=i$. We have $q \dim
C^{\alpha}(G) =(q\dim A)^{b_0([G:s])}(q\dim B)^{b_1([G:s])}$ which
is exactly the contribution of the state $s$ in
$\widehat{T}(G;x,y)$. This proves the equation.
\end{proof}

\section{More general constructions} \label{sec4}

Our construction can be made more general. Let $A$ be any
commutative bigraded ring. Let $B$ be any bigraded module over $\Z$.
For both $A$ and $B$, we assume that the dimension of the space of
homogeneous elements at each degree is finite so that $q\dim A$ and
$q\dim B$ are well-defined as 2--variable power series. For each
integer $k\geq 0$, let $f_k\co  B^{\otimes k} \rightarrow B^{\otimes
k+1}$ be a degree preserving module homomorphism. Given such $A, B$
and $f_k$,  we can construct homology groups in the following
manner.

The chain groups are defined similarly as before.  We fix an
ordering on $E(G)$ and denote the edges by $e_1, \cdots, e_n$. Let
$\alpha$ be a vertex of the cube $\{ 0, 1\}^{E(G)} = \{ 0, 1\}^n$.
Let  $s\subset E(G)$ be the edge set corresponding to $\alpha$. We
assign a copy of $A$ to each component of $[G:s]$ and then take
tensor product over the components. Let $A^{\alpha}(G)$ be the
resulting graded $\mathbb{Z}$--module, with the induced grading
from $A$.   Let $B^{\alpha}(G)=B^{\otimes b_1([G:s])}$. We define
$C^{\alpha}(G)=A^{\alpha}(G)\otimes B^{\alpha}(G)$. Then we define
 $C^i(G):=\oplus_{|\alpha|=i}
C^{\alpha}(G)$.

The differential maps $d^i$ are defined using the multiplication on
$A$ and the homomorphism $f_k$. First, we describe the per-edge map
$d_{\xi}\co  C^{\alpha_1}(G)\rightarrow C^{\alpha_2}(G)$. Let $e$ be
the corresponding edge in $G$ so that $s_2=s_1\cup \{ e\}$. If $e$
joins a component of $[G:s_1]$ to itself, we define $d_{\xi}^A\co 
A^{\alpha_1}(G)\rightarrow A^{\alpha_2}(G)$ to be the identity map,
and define $d_{\xi}^B\co  B^{\alpha_1}(G)\rightarrow B^{\alpha_2}(G)$
to be the homomorphism $f_k\co  B^{\otimes k} \rightarrow B^{\otimes
k+1}$ where $B^{\alpha_1}(G)=B^{\otimes k}$ and
$B^{\alpha_2}(G)=B^{\otimes k+1}$. The per-edge map $d_{\xi}:
C^{\alpha_1}(G)\rightarrow
 C^{\alpha_2}(G)$ is defined to be $d_{\xi}^A \otimes d_{\xi}^B$.
If  $e$ joins two different components, say $E_1$ and $E_2$, of
$[G:s_1]$, we define $d_{\xi}$ to be the multiplication map $A
\otimes A \rightarrow A$ on tensor factors coming from $E_1$ and
$E_2$, and $d_{\xi}$ to be the identity map on the remaining tensor
factors.

As before, we define the differential $d^i\co C^i(G) \rightarrow
C^{i+1}(G)$ by $d^i = \sum_{|\xi|=i} (-1)^{\xi} d_{\xi}$, where
$(-1)^{\xi} =(-1)^{\sum_{i<j} \xi_i} $ and $j$ is the position of
$*$ in $\xi$.

A similar argument as before proves:

\begin{theorem}\label{more general chain complex}$\phantom{99}$
 
{\rm(a)}\qua $0 \rightarrow C^0(G)
\stackrel{d^0}{\rightarrow}C^1(G) \stackrel{d^1}{\rightarrow} \cdots
\stackrel{d^{n-1}}{\rightarrow} C^n(G) \rightarrow 0$ is a chain
complex of bigraded modules whose differential
is degree preserving. Denote this chain complex by $C(G)=C_{A, B, f_k}(G)$.

{\rm(b)}\qua The cohomology groups $H^i(G) (=H^i_{A, B, f_k}(G))$ are
independent of the ordering of
the edges of $G$, and therefore are invariants of the graph $G$.
In fact, the isomorphism type of the graded chain complex $C(G)$ is an
invariant of $G$.

{\rm(c)}\qua The graded Euler characteristic
\begin{align*}
\chi_{q}(C(G))=&\mathrel{\mathop{\sum }\limits_{0\leq i\leq
n}}(-1)^{i} q\dim (H^{i})=\mathrel{\mathop{\sum }\limits_{0\leq
i\leq n}} (-1)^{i} q\dim (C^{i})\\
=&\widehat{T}(G;q\dim A -1,q\dim B
-1).\end{align*}
\end{theorem}

\begin{remark} Some special choices of $A$, $B$ and $f_k$ are as follows.
\vspace{-2pt}

(a)\qua Given any $A$ and $B$ satisfying the basic conditions as
above. Let $b_0$ be a fixed element in $B$ with $\deg b_0 =
(0,0)$. Then we can define $f_k\co  B^{\otimes k} \rightarrow
B^{\otimes (k+1)}$ by $f_k(b)=b\otimes b_0$ for all $b\in
B^{\otimes k}$. In particular, if $A$ and $B$ are given as in the
previous section, and $b_0=1$,
we obtain the construction in the previous section.
\vspace{-2pt}

(b)\qua Let $A, B, b_0$ be as in (a) above. A specific choice of $b_0$
 is $b_0=0$. 
\vspace{-2pt}

(c)\qua Let $B=\mathbb{Z}$, with $q\dim B=1$ and $f_k(b)=b \otimes 1$.
Then $\chi_{q}(C(G))=P_G(q\dim
 A)$. The homology groups are isomorphic to the ones in \cite{HR05}.
\end{remark}
\vspace{-2pt}

\section{Exact sequences}\label{sec5}
\vspace{-2pt}

In this section, we show that our homology groups satisfy a long
exact sequence which can be considered as a categorification for the
deletion-contraction rule. Since we will work with $\widehat{T}(G)$
rather than $T(G)$, we need the deletion contraction rule for
$\widehat{T}(G)$, which we establish in \fullref{sec5.1}. This naturally
leads us to two cases: adding an non-loop edge or adding a loop. The
exact sequences for these two cases are discussed in \fullref{sec5.2} and
\fullref{sec5.3}. respectively.
\vspace{-2pt}

\subsection[Deletion-contraction rule for T-hat(G)]{Deletion-contraction rule for $\widehat{T}(G)$}\label{sec5.1}
\vspace{-2pt}
Let $e$ be a fixed edge of $G$. We wish to understand relations
between $\widehat{T}(G), \widehat{T}(G-e)$, and $\widehat{T}(G/e)$.
Recall that $\widehat{T}(G;x,y)=\sum_{s\subset E(G)} (-1)^{|s|}
(1+x)^{b_0([G:s])} (1+y)^{b_1([G:s])}$.  We have
$\widehat{T}(G;x,y)=\sum_1 + \sum_2$, where $\sum_1$ consists of
terms with those $s$ that do not contain $e$, and $\sum_2$ consists
of the remaining terms.
\vspace{-2pt}

For the first summation, we have $$\textstyle \sum_1 =\sum_{e\not \in s\subset
E(G)} (-1)^{|s|} (1+x)^{b_0([G:s])} (1+y)^{b_1([G:s])}.$$  This
summation is the same as the summation over all $s\subset E(G-e)$.
Furthermore, for each such $s$, $[G:s]=[G-e:s]$. It follows that
$\sum_1 = \widehat{T}(G-e;x,y)$. Note that this is true for all
choices of $e$, including loops and isthmuses.
\vspace{-2pt}

For the second summation, we have $\sum_2 =\sum_{e \in s\subset
E(G)} (-1)^{|s|} (1+x)^{b_0([G:s])} (1+y)^{b_1([G:s])}$. Each such
$s$ can be written as $s=s_1\cup \{ e\}$, where $s_1$ corresponds to
a subset of $E(G/e)$, which we also denote by $s_1$. Note that
$(-1)^{|s|}=-(-1)^{|s_1|}$ since $|s|=|s_1|+1$. If we assume that
$e$ is not a loop, we have $[G:s]\simeq [G/e:s_1]$ where $\simeq$
stands for homotopically equivalent. It follows that
$\sum_2=-\widehat{T}(G/e, x,y)$. Note that this is true for all $e$
that are not loops (therefore $e$ can be an isthmus).
\vspace{-2pt}

We have proved (1) of the following:

\begin{proposition} \label{d-c rule for T hat}
Let $e$ be an edge in a graph $G$.

{\rm(1)}\qua If $e$ is not a loop (but possibly an isthmus), then
$\widehat{T}(G;x,y)=\widehat{T}(G-e;x,y)-\widehat{T}(G/e;x,y)$.

{\rm(2)}\qua If $e$ is a loop, then
$\widehat{T}(G;x,y)=\widehat{T}(G-e;x,y)-(1+y)\widehat{T}(G/e;x,y)=-y
\widehat{T}(G/e;x,y)$.
Of course, we have $G-e=G/e$.

{\rm(3)}\qua If $e$ is an isthmus, then
$$\widehat{T}(G;x,y)=\widehat{T}(G-e;x,y)-\widehat{T}(G/e;x,y)
=x\widehat{T}(G/e;x,y).$$
We also have $\widehat{T}(G-e;x,y)=(1+x)\widehat{T}(G/e;x,y)$.
\end{proposition}

\begin{proof} Part (1) is proved above. For (2), we claim
$\sum_2= - (1+y) \widehat{T}(G/e, x,y)$. This is because $[G:s]$ is
obtained from $[G/e:s_1]$ by adding the loop $e$ to an existing
vertex, which then implies $b_1([G:s])=b_1([G/e:s_1])+1$. For (3),
We need only to show
$\widehat{T}(G-e;x,y)=(1+x)\widehat{T}(G/e;x,y)$. Both sides can be
written as a summation over $s\in E(G-e)=E(G/e)$. The only
difference is that $b_0([G-e:s])=b_0([G/e:s])+1$. This implies our
equation.
\end{proof}

\subsection{Exact sequence for non-loop edges}\label{sec5.2}

Let $e$ be an edge of a graph $G$. Assume that $e$ is not a loop. We
will construct degree preserving chain maps $\alpha:
C^{i-1}(G/e)\rightarrow C^i(G)$, and $\beta\co  C^i(G)\rightarrow
C^i(G-e)$ such that
\[ 0\rightarrow C^{i-1}(G/e) \overset{\alpha}{\rightarrow} C^i(G)
\overset{\beta}{\rightarrow} C^{i}(G-e)\rightarrow 0 \] is exact.

First, we describe $\alpha$. It is enough to define
$\alpha|_{C^s(G/e)}$ for all $s\subset E(G/e)$ with $|s|=i-1$. Let
$s\subset E(G/e)$ with $|s|=i-1$. Let $s_e=s\cup\{ e\}$. Then
$s_e\subset E(G)$ with $|s_e|=i$. Since $e$ is not a loop, the two
graphs $[G:s_e]$ and $[G/e:s]$ are homotopically equivalent under
the natural map that contracts $e$ to a point. This induces a
natural isomorphism from $C^{s}(G/e)$ to $C^{s_e}(G)$. Define our
map $\alpha|_{C^s(G)}\co  C^{s}(G/e) \rightarrow C^{s_e}(G) $ to be
this isomorphism. Taking summation over $s$, we obtain the map
$\alpha\co  C^{i-1}(G/e) \rightarrow C^i(G)$.

Next, we describe the map $\beta$. Just as before, it is enough to
define $\beta|_{C^s(G)}$ for all $s\subset E(G)$ with $|s|=i$. We
consider two cases. If $e \in s$, we define $\beta|_{C^s(G)}$ to be
the zero map. If $e\not \in s$, the two graphs $[G-e:s]$ and $[G:s]$
are identical, and therefore the groups $C^s(G-e)$ and $C^s(G)$ are
naturally identified. We define $\beta|_{C^s(G)}$ to be this
identity isomorphism composed with the inclusion map.

A diagram chasing argument shows:

\begin{theorem} \label{exact seq. no-loop}$\phantom{99}$

{\rm(a)}\qua If $e$ is not a loop, the above
defines an exact sequence of chain maps, and therefore

{\rm(b)}\qua it induces a long exact sequence:
\begin{align*}
0\rightarrow&
H^{0}(G)\overset{\beta ^{\ast }}{\rightarrow }H^{0}(G-e)
\overset{\gamma ^{\ast }}{\rightarrow }H^{0}(G/e)
\overset{\alpha ^{\ast }}{ \rightarrow }\\
&H^{1}(G)\overset{\beta
^{\ast }}{\rightarrow }H^{1}(G-e) \overset{\gamma ^{\ast
}}{\rightarrow }H^{1}(G/e)
\overset{\alpha ^{\ast }}{ \rightarrow }\ldots 
\end{align*}
\end{theorem}

It will be useful for further computations to understand the action
of the map $\gamma ^*$. This follows from tracing back the elements
involved in the diagram chasing of the zig-zag lemma. The result is
as follows.

\begin{remark}\label{interpretation gamma} The connecting homomorphism $\gamma^*\co  H^i(G-e) \rightarrow
H^i(G/e)$ acts as follows. Each cycle $z \in C^i(G-e)$ is a linear
combination of states for the graph $G-e$. That is, $z=\sum n_k
(s_k, c_k)$ where  $n_k \in \Z, s_k$ is a subset of $E(G-e)$, $c_k
\in C^{s_k}(G-e)$. We add the edge $e$ to each $s_k$ to get $s_k\cup
\{ e\}$, and replace $c_k$ by $d_{\xi}(c_k)$ where $d_{\xi}$ is
described in \fullref{sec3.2}. More specifically, $d_{\xi}(c_k)$ is
obtained using the multiplication on $A$ if $e$ joins two components
of $G-e$, otherwise $d_{\xi}(c_k)$ is obtained using the map
$B\rightarrow B\otimes B$ sending $b$ to $b\otimes 1$.
\end{remark}

\subsection{Exact sequence for loops}\label{sec5.3}
Now we assume that $e$ is a loop. We will define a new chain complex
denoted by $C(G/e)\otimes B$. Then we will show that there is a
short exact sequence of chain maps:
\[ 0\rightarrow C^{i-1}(G/e)\otimes B \overset{\alpha}{\rightarrow} C^i(G)
\overset{\beta}{\rightarrow} C^{i}(G-e)\rightarrow 0 \]
Of course, $G-e=G/e$.

First, we describe $C(G/e)\otimes B$. For each $i$, define its
$i^{th}$ chain group to be $C^i(G/e)\otimes B$. Next, define
differential $d\co  C^{i-1}(G/e)\otimes B \rightarrow C^i(G/e)\otimes
B$ in terms of per-edge maps. Let $s\subset E(G/e)$, with $|s|=i-1$,
and $f\in E(G/e)-s$. The per-edge map $d\co  C^s(G/e)\otimes B
\rightarrow C^{s\cup \{ f\}}(G/e)\otimes B$ is defined as follows.
If $f$ connects two different components of $[G/e:s]$, $d$ is
defined the same way as before by multiplying the ``coloring" on
these two components. If $f$ joins a component of $[G/e:s]$ to
itself, $C^{s\cup \{ f\}}(G/e)\otimes B\cong (C^s(G/e)\otimes B)
\otimes B$ via a natural isomorphism, and then per-edge map $d$ is
defined to be $d(x)=x\otimes 1$ for each $x\in C^s(G/e)\otimes B$.
Note that this differential map is different from just taking tensor
product with the identity, that is: $d_{C(G/e)\otimes B} \ne
d_{C(G/e)}\otimes Id_B$.

We now describe the maps $\alpha$ and $\beta$. The map $\alpha:
C^{i-1}(G/e)\otimes B \rightarrow C^i(G)$ is defined as follows. For
each $s\subset E(G/e)$ with $i$ edges, $s\cup \{ e\} \in E(G)$ has
$i+1$ edges. The graph $[G:s\cup \{ e\}]$ is obtained from $[G:s]$
by adding the loop $e$. Thus $C^s(G/e)\otimes B $ is naturally
isomorphic to $C^{s\cup \{ e\}}(G)$. Define $\alpha$ to be this
isomorphism summed over all $s$. The map $\beta: C^i(G)\rightarrow
C^i(G-e)$ is defined the same way as in the above section. In other
words, it is the projection map that kills all summands $C^s(G)$
where $e\in s$.

\begin{theorem}
Let $e$ be a loop in $G$. Then:

{\rm(a)}\qua There is a short exact sequence of chain maps
\[ 0\rightarrow C^{i-1}(G/e)\otimes B \overset{\alpha}{\rightarrow} C^i(G)
\overset{\beta}{\rightarrow} C^{i}(G-e)\rightarrow 0 \] described
above. Therefore:

{\rm(b)}\qua This induces a long exact sequence: 
\begin{align*}
0\rightarrow&
H^{0}(G)\overset{\beta ^{\ast }}{\rightarrow }H^{0}(G-e)
\overset{\gamma ^{\ast }}{\rightarrow }H^{0}(C(G/e)\otimes
B)\overset{\alpha ^{\ast }}{ \rightarrow }\\
&H^{1}(G)\overset{\beta
^{\ast }}{\rightarrow }H^{1}(G-e) \overset{\gamma ^{\ast
}}{\rightarrow }H^{1}(C(G/e)\otimes B)
\overset{\alpha ^{\ast }}{ \rightarrow }\ldots 
\end{align*}
\end{theorem}

The proof consists of a standard but tedious diagram chasing
argument, and is left as an exercise.

\section{Other properties}\label{sec6}

We prove some other properties of our cohomology groups.

\subsection{Adding a pendant edge}
As an application of the above long exact sequence, let us consider
the effect of adding a pendant edge on the cohomology groups. Recall
that a {\em pendant vertex} in a graph is a vertex of degree one,
and a {\em pendant edge} is an edge connecting a pendant vertex to
another vertex. Let $e$ be a pendant edge in a graph $G$, then
\fullref{d-c rule for T hat} implies
$\widehat{T}(G;x,y)=x\widehat{T}(G/e;x,y)$. On the level of our
cohomology groups, we have the following:

\begin{theorem}\label{pendant edge} If $e$ is a pendant edge of
the graph $G$, then $H^k(G)\cong H^k(G/e)\{ (1,0)\}$,  where $\{
(1,0)\} $ is the operation that shifts the degree up by $(1,0)$.
\end{theorem}

\begin{remark}
Let $A'=\mathbb{Z}x$ be the submodule generated by $x$. Then $A=\Z 1
\oplus \Z x$, and $A'\cong \Z \{(1,0)\}$ as bigraded modules. The
above can be rephrased as
\[ H^k(G)\cong H^k(G/e)\otimes A'\]
This equation works for more general algebras $A$ satisfying $A=\Z 1
\oplus A'$.
\end{remark}

\begin{proof}[Proof of \fullref{pendant edge}] The proof is similar to the proof of
an analogous result for the chromatic cohomology (Theorem 22 in
\cite{HR04}).

Consider the operations of contracting and deleting $e$ in $G$.
Denote the graph $G/e$ by $G_1$. We have $G/e=G_1,$ and
$G-e=G_1\sqcup \{v\}$, where $v$ is the end point of $e$ with $\deg
v = 1$. Consider the exact sequence
\begin{equation*}
\cdots \rightarrow H^{i-1}(G_1\sqcup
\{v\})\overset{\gamma^*}{\rightarrow } H^{i-1}(G_1)
 \overset{\alpha^*}{\rightarrow } H^{i}(G)\overset{\beta^*}{\rightarrow
} H^{i}(G_1 \sqcup \{v\})\overset{\gamma^*}{\rightarrow }
H^{i}(G_1)\rightarrow \cdots
\end{equation*}
We need to understand the map
\begin{equation*}
H^{i}(G_1\sqcup \{v\})\overset{\gamma^*}{\rightarrow }H^{i}(G_1)
\end{equation*}
It is easy to understand the impact of adding an isolated vertex on
the cohomology groups. We have
\begin{equation*}
H^{i}(G_{1}\sqcup \{v\})\cong H^{i}(G_{1})\otimes A.
\end{equation*}
This can be seen by first noting the same equation holds on the
level of chain groups. Furthermore,  the differential map restricted
on the tensor factor $A$ is the identity map. This implies that
equation holds on the level of homology groups.

We therefore identify $H^{i}(G_{1}\sqcup \{v\})$ with
$H^{i}(G_{1})\otimes A$. The map $\gamma^{\ast} \co  H^{i}(G_1\sqcup
\{v\})\rightarrow H^{i}(G_1)$ sends $u\otimes 1$ to $(-1)^i u$. In
particular, $\gamma^*$ is onto. Therefore, the above long exact
sequence becomes a collection of short exact sequences:
\begin{equation}
0\rightarrow H^{i}(G)\overset{\beta^*}{\rightarrow }
H^{i}(G_{1}\sqcup \{v\})\overset{\gamma^*}{\rightarrow }
H^{i}(G_{1})\rightarrow 0
\end{equation}
Hence, $H^i(G)\cong \ker \gamma^* $. We define a homomorphism:
\[ f\co  H^i(G_1)\otimes A' \rightarrow \ker \gamma^* \mbox {  by } \\
 f(u\otimes a')=u\otimes a' -(-1)^i \gamma^*(u\otimes a') \otimes 1 \]
One checks that $f$ is an isomorphism of $\Z$--modules. Therefore,
$H^i(G)\cong \ker \gamma^* \cong H^i(G_1)\otimes A'$.
\end{proof}

\begin{remark} \label{generators pendant}
The above gives a complete description of the generators of
$H^i(G)$. For the specific case of $A'= \Z x$, $H^i(G)$ is spanned
by $ \alpha \otimes x -(-1)^i \gamma^*(\alpha \otimes x) \otimes 1$,
where $\alpha$ ranges over the generators of $H^i(G_1)$.
\end{remark}
\newpage
Using inductively \fullref{pendant edge} on the number of edges
implies:
\vspace{-6pt}

\begin{corollary}\label{tree}
If $G=T_n$ is a tree with $n$ edges, then $H^0(G)\cong \mathbb{Z}\{
(n,0)\} \oplus \Z \{ (n+1,0)\}$, and $H^i(G)=0$ for all $i>0$.
\end{corollary}
\vspace{-6pt}

We remark that, for the more general construction in \fullref{sec4}, a statement similar to \fullref{tree}  still holds with $A'$ being the direct summand of $A$
satisfying $A= \Z 1 \oplus A'$
\vspace{-6pt}

\subsection{Overlapping with the chromatic cohomology}
\vspace{-6pt}

\fullref{tree} shows that, for all trees, the Tutte cohomology
and the chromatic cohomology agree. This is extended to general
graphs in part (a) below.
\vspace{-6pt}

\begin{theorem} Let $\ell$ be the length of the shortest cycle in a graph $G$.

{\rm(a)}\qua For all $i<
 \ell -1$, we have
 \[ H^i_{\mbox{\footnotesize\rm Tutte}}(G)\cong H^i_{\mbox{\footnotesize\rm chromatic}}(G)\]
{\rm(b)}\qua If $G$ is the Tait graph of an unoriented framed alternating
link diagram $D$, then for all $i<\ell -1$,
\[ H^{i,j,0}_{\mbox{\footnotesize\rm Tutte}}(G)\cong H_{p,q}(D) \]
where $p=|V(G)|-i-2j, q=|E(G)|-2|V(G)|+4j,$ and $H_{p,q}$ is the
homology groups of the version of Khovanov cohomology theory for
unoriented framed link defined by Viro in \cite{V04}.
\end{theorem}
\vspace{-6pt}

\begin{proof}
(a)\qua For each such $i$, the algebra $B$ is not involved in the chain
groups $C^{0}_{\mbox{\footnotesize Tutte}}(G), \cdots,
C^{i+1}_{\mbox{\footnotesize Tutte}}(G)$ since subgraphs involved,
$[G:s]$, contains no cycles. It follows that the chain groups and
differentials in the Tutte cohomology and the chain groups in the
chromatic cohomology are all the same. Therefore the corresponding
cohomology groups are isomorphic.
\vspace{-6pt}

(b)\qua This is a consequence of
Theorem 24 of Helme-Guizon, Przytycki and Rong \cite{HPR} where a
similar relation between chromatic cohomology and Khovanov's link
homology was established.
\end{proof}
\vspace{-6pt}

\subsection{Functorial property}
The classical homology theory is a functor: continuous maps between
spaces induce homomorphism between homology groups. Khovanov's link
homology $Kh$ also satisfies a functorial property: a cobordism $C$
between two links $L_1$ and $L_2$ induces a homomorphism $Kh(C)\co 
Kh(L_1) \rightarrow Kh(L_2)$, well defined up to sign
(Jacobsson \cite{J02}, Khovanov \cite{K02}).

 For our Tutte cohomology, we can associate homomorphisms between
 homology groups to each inclusion map of graphs. Essentially, it
 is the iteration of the map $\beta^*$ in the above long exact
 sequence, with some additional attention to any possible
 additional vertices in the ambient graph. More specifically,
let $K$ be a subgraph of $G$, and we denote this relation by
$K\subseteq G$. Define $\beta\co C^i(G) \rightarrow C^i(K)$ as
follows. Since $C^i(G):=\oplus_{|s|=i} C^s(G)$, it is enough to
define $f|C^s(G)$.

(1)\qua If $s\not \subseteq E(K)$, we define $\beta|C^s(G)=0$.

(2)\qua If $s\subseteq E(K)$, then $[G:s]$ is $[K:s]$ union $l$ vertices
of $G$ where $l\geq 0$. We have $C^s(G)\cong C^s(K)\otimes
A^{\otimes l}$. Since $A = \Z 1\oplus \Z x$, $A^{\otimes l}=A_0
\oplus A_1$, where $A_0\cong \Z$ is the direct summand generated by
$1 \otimes \cdots \otimes 1$. This implies $C^s(G)\cong C^s(K)
\otimes A_0 \oplus C^s(K) \otimes A_1$. We define $\beta(g\otimes
1\otimes \cdots \otimes 1)=g, \beta(g\otimes a)=0$ for all $g\in
C^s(K), a \in A_1$. 

\begin{theorem}$\phantom{99}$

{\rm(1)}\qua $\beta\co  C(G) \rightarrow C(K)$ is a degree preserving
chain map, therefore

{\rm(2)\qua} it induces a degree preserving homomorphism $\beta^*\co  H^i(G)
\rightarrow H^i(K)$.

{\rm(3)}\qua This correspondence is natural. That is,
if $L\subseteq K, K\subseteq G$ are subgraphs, then the diagram
$$
\xymatrix{
H^*(G) \ar[rr]^{\beta^*} \ar[dr]_{\beta^*_2}
 & &H^*(L) \\
& H^*(K) \ar[ur]_{\beta^*_1} & }
$$
commutes, where $\beta^*, \beta^*_1, \beta^*_2 $ are the
homomorphisms induced by the inclusions $L\subseteq G, L\subseteq
K,$ and $K\subseteq G$ respectively.
\end{theorem}

\section{Examples}\label{sec7}

\begin{exa}\label{ex1}
Let $L_1$ be the graph with one vertex and one
loop, that is, $L_1 = $ \psfig{figure=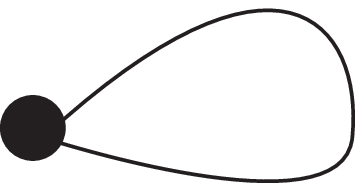,height=.3 cm}.
Our construction yields the chain complex:
$$
0 \rightarrow A \overset{d^0}{\rightarrow} A \otimes B {\rightarrow}
0
$$
The differential $d^0$ maps $1 \mapsto 1_A\otimes 1_B$ and $x
\mapsto x \otimes 1_B$. It is easy to see that
$$
 H^1(L_1)= \langle 1_A \otimes y, x\otimes y \rangle \cong A\otimes B'\cong \Z \{(0,1)\} \oplus
\mathbb{Z}\{ (1,1)\} , H^i(L_1)=0 \mbox{ for $i\neq 1.$ }
$$
where $B'=\Z \{ y\}$ satisfies $B\cong \mbox{ span}\{ 1\} \oplus
B'$. Hence we have $\chi(H^*(L_1))= -y-xy = \widehat{T}(L_1;x,y)$
\end{exa}

\begin{exa}\label{ex2} Let $P_2 = $
\psfig{figure=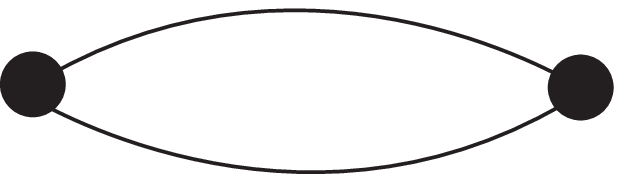,height=.3 cm},
 that is, the Polygon with two sides. The corresponding chain complex is:
$$
0 \rightarrow A \otimes A \overset{d^0}{\rightarrow} A \oplus A
\overset{d^1}{\rightarrow} A \otimes B \overset{d^2}{\rightarrow} 0
$$
Where $d^0$ maps $x\otimes x \mapsto (0, 0)$, $1_A \otimes 1_A
\mapsto (1_A, 1_A)$, $x \otimes 1_A$ and $1_A \otimes x \mapsto
(x,x)$. The kernel of $d^0$ is generated by the elements $x
\otimes x$ and $x \otimes 1_A- 1_A \otimes x$. Thus $H^0 (P_2)
\cong A\{ 1\} \cong \Z(1,0)\oplus \Z(2,0)$. We also have $d^1$
maps $(1_A,0)\mapsto -1_A \otimes 1_B$, $(x,0) \mapsto -x \otimes
1_B$, $(0,1_A) \mapsto 1_A \otimes 1_B $ and $ (0,x) \mapsto  x
\otimes 1_B $. This implies $H^1(P_2)=0$ and $H^2(P_2) \cong
A\otimes B'\cong \Z\{(0,1)\}\oplus \Z \{(1,1)\}$. Clearly
$H^i(P_2)=0$, for $i\geq 3$. Hence
$\chi(H^*(P_2))= x^2+x+y+xy = \widehat{T}(P_2;x,y)$.
\end{exa}

\begin{exa}\label{ex3} Let $L_2=$ \psfig{figure=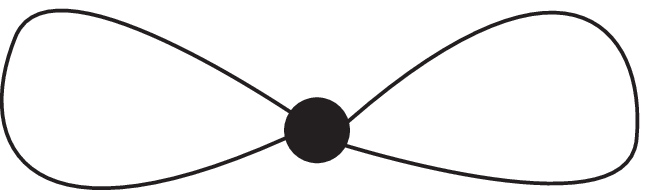,height=.3
cm} (ie one vertex and two loops). Its chain complex is:
$$
0 \rightarrow A \overset{d^0}{\rightarrow} A \otimes B \oplus A
\otimes B \overset{d^1}{\rightarrow} A \otimes B \otimes B
\rightarrow 0
$$
Where $d^0$ maps $1 \mapsto (1_A\otimes 1_B,1_A\otimes 1_B)$ and $x
\mapsto (x \otimes 1_B, x \otimes 1_B)$, $d_1$ maps $(a \otimes b,
0) \mapsto - a \otimes b \otimes 1_B$ and $(0,a \otimes b) \mapsto
a \otimes b \otimes 1_B$. Hence, $H^1(L_2)=\langle (1_A \otimes y,
1_A \otimes y)\rangle \oplus \langle(x \otimes y, x \otimes
y)\rangle \cong A \otimes B' \cong \Z \{(0,1)\} \oplus \Z
\{(1,1)\}$, $H^2 \cong A \otimes B \otimes B' \cong \Z \{(0,1)\}
\oplus \Z \{(1,1)\}\oplus \Z \{(0,2)\} \oplus \Z
\{(1,2)\}$ and $H^i (L_2)=0$ for $i \neq 1,2$.

Clearly $\chi(H^*(L_2))= xy^2+y^2 = \widehat{T}(L_2;x,y)$.
\end{exa}

\begin{exa}\label{ex4} Let $G=$ \psfig{figure=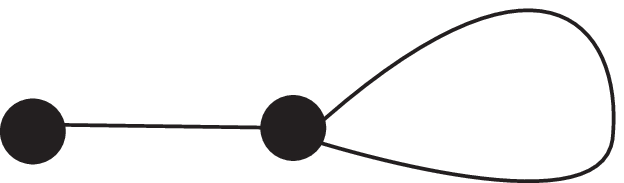,height=.3
cm}, obtained from $L_1$ by adding a pendant edge. The
computations in \fullref{ex1} and \fullref{generators pendant}
yield:
\begin{eqnarray*} \label{gen hi(l2)}
H^1(G)&=& \langle1_A\otimes x \otimes y - x\otimes 1_A \otimes
y\rangle \oplus \langle x\otimes x \otimes y\rangle\\
&&\qquad\qquad\qquad \cong A\otimes
B'\{(1,0)\}\cong \mathbb{Z}\{(1,1)\} \oplus \mathbb{Z}\{ (2,1)\} ,\\
H^i(G)&=&0 \mbox{\ \ for\ \ $i\neq 1.$ } 
\end{eqnarray*}
We have that $\chi(H^*(G))= -xy-x^2y = \widehat{T}(G;x,y)$.
\end{exa}

\begin{exa}\label{ex5}  Let $G = $
\psfig{figure=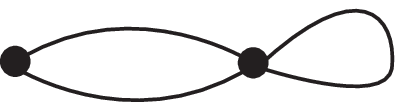,height=.3 cm}, whose construction of
chain complex and differential was described in \fullref{sec3.2}. An
easy way to compute the cohomology groups for this graph is using
the exact sequence given in \fullref{exact seq. no-loop}.  We choose the edge $e$ to be one
\labellist\tiny\hair1pt
\pinlabel $e$ [b] at 45 25
\endlabellist 
of the parallel edges on
$G:$\psfig{figure=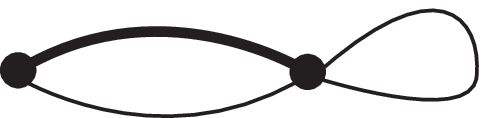,height=.5 cm}. Then $G-e$ is the
graph in \fullref{ex4}, and $G/e$ is $L_2$ as in \fullref{ex3}.

According to \fullref{exact seq. no-loop} we have the following long
exact sequence:
$$0\rightarrow
H^{0}(G)\overset{\beta ^{\ast }}{\rightarrow }H^{0}(G-e)
\overset{\gamma ^{\ast }}{\rightarrow }H^{0}(L_2)\overset{\alpha
^{\ast }}{ \rightarrow }H^{1}(G)\overset{\beta ^{\ast }}{\rightarrow
}H^{1}(G-e) \overset{\gamma ^{\ast }}{\rightarrow
}H^{1}(L_2)\overset{\alpha ^{\ast }}{\rightarrow} \ldots$$
If we use the computations in Examples \ref{ex3} and \ref{ex4}, then the long
exact sequence above decomposes into the following exact
sequences:
\begin{gather*}
0\rightarrow H^{0}(G)\overset{}{\rightarrow } 0
\\
0 { \rightarrow }H^{1}(G)\overset{\beta ^{\ast }}{\rightarrow }
A\otimes B'\{(1,0)\} \overset{\gamma ^{\ast }}{\rightarrow }
A\otimes B' \overset{\alpha ^{\ast }}{\rightarrow}
H^2(G)\overset{\beta ^{\ast }}{\rightarrow }0
\\
0\rightarrow A \otimes B \otimes B' \overset{\alpha ^*}{\rightarrow}
H^{3}(G)\overset{}{\rightarrow } 0
\end{gather*}
Hence we have $H^0(G) \cong 0$, $H^3(G) \cong A \otimes B \otimes B'
\cong \Z\{(0,1)\} \otimes \Z \{ (1,1) \} \otimes \Z \{ (0,2) \}
\otimes \Z \{ (1,2) \}$. To compute $H^i(G)$ for $i=1,2$, we need to
understand the map $\gamma^*\co  H^1(G-e) \rightarrow H^1(L_2)$. First
note that, by \fullref{ex4}, $H^1(G-e)\cong A \otimes B'\{(0,1)\}$ has
basis $\{1\otimes x \otimes y - x\otimes 1 \otimes y, x\otimes x
\otimes y \}$ and, by \fullref{ex3}, $H^1(L_2) \cong A \otimes B'$ has
basis $\{ ( 1\otimes y, 1 \otimes y ), ( x\otimes y, x \otimes y
)\}$. By \fullref{interpretation gamma} we have $\gamma^*(
1\otimes x \otimes y - x\otimes 1 \otimes y)= x\otimes y - x\otimes
y =0$ and $\gamma^*( x \otimes x \otimes y)= x^2 \otimes y=0$
(recall $x^2=0$ in $A$). Thus $\gamma^*$ is the zero map. Hence
$H^{1}(G) \cong H^1(G-e) \cong \Z\{(2,1)\} \oplus \Z\{(1,1)\}$ and
$H^2(G) \cong A\otimes B' \cong \Z\{(0,1)\} \oplus \Z\{(1,1)\}$.

The graded Euler characteristic for the cohomology groups for $G$ is
$-(1+x)(1+y)y+(1+x)y-(1+x)xy =-y(1+x)(x+y)$, which agrees with
$\widehat{T}(G;x,y)$.
\end{exa}

\begin{exa}\label{ex6} Let $K_3$ be the complete graph with 3 vertices.
To compute its cohomology groups we will use the exact sequence
given in \ref{exact seq. no-loop}. In this case $G-e= T_2$, the
tree with two edges and $G/e=P_2$ as in \fullref{ex2} above. We get:
$$0\rightarrow
H^{0}(G)\overset{\beta ^{\ast }}{\rightarrow }H^{0}(T_2)
\overset{\gamma ^{\ast }}{\rightarrow }H^{0}(P_2)\overset{\alpha
^{\ast }}{ \rightarrow }H^{1}(G)\overset{\beta ^{\ast }}{\rightarrow
}H^{1}(T_2) \overset{\gamma ^{\ast }}{\rightarrow
}H^{1}(P_2)\rightarrow \ldots $$ Using the information of the
previous examples, if we substitute the homology groups that are
zero, we end up with the following exact sequences:
\begin{gather*}
0\rightarrow H^{0}(G)\overset{\beta ^{\ast }}{\rightarrow
}H^{0}(T_2) \overset{\gamma ^{\ast }}{\rightarrow
}H^{0}(P_2)\overset{\alpha ^{\ast }}{ \rightarrow
}H^{1}(G)\rightarrow 0 
\\
0\rightarrow H^2(P_2)\overset{\alpha^{\ast }}{ \rightarrow
}H^{3}(G)\rightarrow 0
\end{gather*}
As a consequence, we have that $H^0 (G) \cong \ker \gamma^{\ast}$,
$H^1 (G) \cong H^0 (P_2) / {\rm Im} \gamma ^{\ast}$, $H^2 (G) \cong 0$ and
$H^3 (G) \cong H^2 (P_2) \cong \Z(0,1) \oplus \Z(1,1)$. So we need
to understand the map $\gamma^\ast\co  H^0 (P_2)\rightarrow H^0 (T_2)$.
The generators for $H^0 (P_2)$ are the elements $1\otimes x\otimes x
- x \otimes 1 \otimes x + x \otimes x \otimes 1$ and $x \otimes x
\otimes x $ which are mapped to $ 2 x\otimes x $ and 0,
respectively. Hence $H^0(G) \cong \Z(3,0)$ and $H^1 \cong \langle
x\otimes 1 - 1 \otimes x \rangle \oplus \langle x \otimes x \rangle
/ \langle2 x \otimes x \rangle \cong \Z(1,0)\oplus \Z_2 (2,0)$ The
graded Euler characteristic is $x^3-x-xy-y=\widehat{T}(G;x,y)$
\end{exa}

\bibliographystyle{gtart}
\bibliography{link}

\begin{thebibliography}{}
\providecommand\bibmarginpar{\leavevmode\marginpar}
\def\urlstyle#1{{\tt #1}}

\bibitem{HPR}
\textbf{L Helme-Guizon}, \textbf{J\,H Przytycki}, \textbf{Y Rong},
  \emph{Torsion in graph homology}, Fund. Math. 190 (2006) 139--177
  \xox{MR}{2232857}

\bibitem{HR05}
\textbf{L Helme-Guizon}, \textbf{Y Rong}, \emph{Graph cohomologies from
  arbitrary algebras} \xox{arXiv}{math.QA/0506023}

\bibitem{HR04}
\textbf{L Helme-Guizon}, \textbf{Y Rong},
  \href{http://dx.doi.org/10.2140/agt.2005.5.1365} {\emph{A categorification
  for the chromatic polynomial}}, Algebr. Geom. Topol. 5 (2005) 1365--1388
  \xox{MR}{2171813}

\bibitem{J02}
\textbf{M Jacobsson}, \href{http://dx.doi.org/10.2140/agt.2004.4.1211}
  {\emph{An invariant of link cobordisms from {K}hovanov homology}}, Algebr.
  Geom. Topol. 4 (2004) 1211--1251 \xox{MR}{2113903}

\bibitem{K00}
\textbf{M Khovanov},
  \href{http://projecteuclid.org/getRecord?id=euclid.dmj/1092749199} {\emph{A
  categorification of the {J}ones polynomial}}, Duke Math. J. 101 (2000)
  359--426 \xox{MR}{1740682}

\bibitem{K02}
\textbf{M Khovanov}, \href{http://dx.doi.org/10.1090/S0002-9947-05-03665-2}
  {\emph{An invariant of tangle cobordisms}}, Trans. Amer. Math. Soc. 358
  (2006) 315--327 \xox{MR}{2171235}

\bibitem{KR04}
\textbf{M Khovanov}, \textbf{L Rozansky}, \emph{Matrix factorizations and link
  homology} \xox{arXiv}{math.QA/0401268}

\bibitem{KR05}
\textbf{M Khovanov}, \textbf{L Rozansky}, \emph{Matrix factorizations and link
  homology II} \xox{arXiv}{math.QA/0505056}

\bibitem{V04}
\textbf{O Viro}, \emph{Khovanov homology, its definitions and ramifications},
  Fund. Math. 184 (2004) 317--342 \xox{MR}{2128056}

\bibitem{W}
\textbf{D\,J\,A Welsh}, \emph{Complexity: knots, colourings and counting},
  London Mathematical Society Lecture Note Series 186, Cambridge University
  Press (1993) \xox{MR}{1245272}

\end{thebibliography}

\end{document}